\documentclass[12pt]{article}
\usepackage{amssymb, amsmath, amsthm, amsfonts}
\usepackage{tikz}
\usetikzlibrary{patterns}
\usetikzlibrary{shapes}
\usepackage{fullpage}
\usepackage{hyperref}
\newtheorem{theorem}{Theorem}
\newtheorem{lemma}[theorem]{Lemma}
\newtheorem{ques}[theorem]{Question}
\newtheorem{claim}[theorem]{Claim}
\newtheorem{prop}[theorem]{Proposition}
\theoremstyle{definition}
\newtheorem{rmk}[theorem]{Remark}
\theoremstyle{definition}
\newtheorem{defn}[theorem]{Definition}
\theoremstyle{definition}
\newtheorem{expl}[theorem]{Example}
\newcommand\ex{\ensuremath{\mathrm{ex}}}

\title{Extremal graphs for blow-ups of cycles and trees}
\author{
Hong Liu
\thanks{Department of Mathematical Sciences,
University of Illinois at Urbana-Champaign, Urbana, Illinois 61801, USA {\tt hliu36@illinois.edu}.}
}
\begin{document}
\maketitle
\begin{abstract}
The \emph{blow-up} of a graph $H$ is the graph obtained from replacing each edge in $H$ by a clique of the same size where the new vertices of the cliques are all different. Erd\H{o}s et al. and Chen et al. determined the extremal number of blow-ups of stars. Glebov determined the extremal number and found all extremal graphs for blow-ups of paths. We determined the extremal number and found the extremal graphs for the blow-ups of cycles and a large class of trees, when $n$ is sufficiently large. This generalizes their results. The additional aim of our note is to draw attention to a powerful tool, a classical decomposition theorem of Simonovits.
\end{abstract}


\section{Introduction}

Notation in this note is standard. We consider undirected graphs without loops and multiedges. For a graph $G$, denote by $E(G)$ the set of edges and $V(G)$ the set of vertices of $G$. The \emph{order} of a graph is the number of its vertices. The number of edges of $G$ is denoted by $e(G) = |E(G)|$. For $U \subseteq V(G)$, let $G[U]$ be the subgraph of $G$ induced by $U$. A path on $k$ vertices is denoted by $P_k$, a star with $k+1$ vertices is denoted by $S_k$ and a cycle with $k$ edges is denoted by $C_k$. A \emph{matching} in $G$ is a set of vertex disjoint edges from $E(G)$, denote by $M_k$ a matching of size $k$. Denote by $T_{n,p}$ the $p$-class \emph{Tur\'an graph}, namely the complete $p$-partite graph on $n$ vertices with the size of each partite set as equal as possible.

The \emph{extremal number}, $\ex(n,H)$, of a graph $H$ is the maximum number of edges in a graph on $n$ vertices which does not contain $H$ as a subgraph. An $H$-free $n$-vertex graph with $\ex(n,H)$ edges is called an \emph{extremal graph} for $H$, or $H$-\emph{extremal}. Tur\'an \cite{Turan41,Turan54} showed that $T_{n,p}$ is the unique extremal graph for $K_{p+1}$. The Erd\H{o}s-Stone-Simonovits Theorem \cite{ErSim66,ErSt46} states that asymptotically Tur\'an's construction is best-possible for any $(p+1)$-chromatic graph $H$ (as long as $p\ge 2$). More precisely $\ex(n,H) = \left(1-\frac{1}{p}\right)\frac{n^2}{2} + o(n^2)$.

\begin{figure}[ht]
\begin{center}
\includegraphics{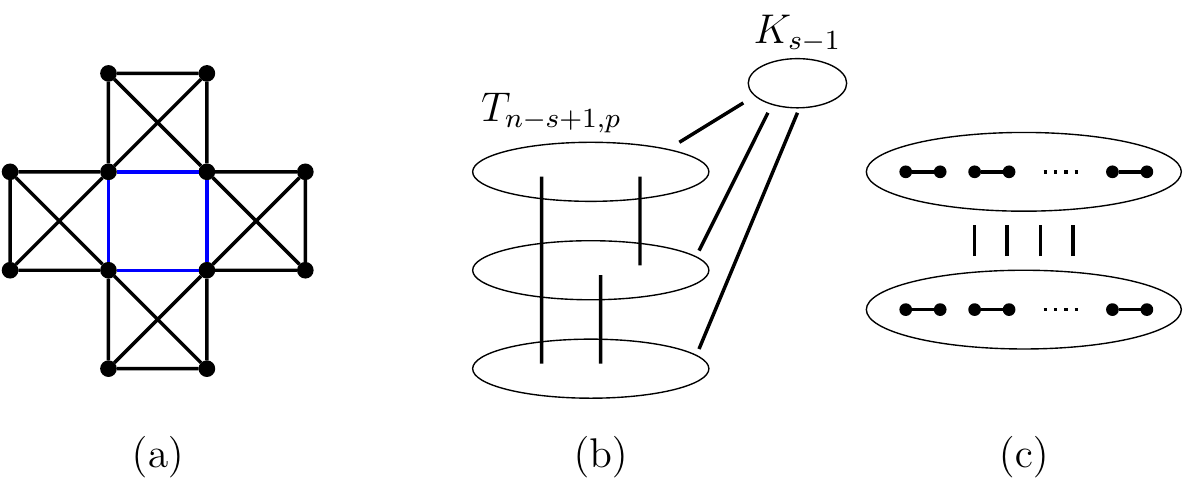}
\end{center}
\caption{(a) $C_4^4$; (b) $H(n,p,s)$; (c) $H^*(n)$.}
\label{fig-0}
\end{figure}

Given a graph $H$, the \emph{blow-up} of $H$, denoted as $H^{p+1}$, is the graph obtained from replacing each edge in $H$ by a clique of size $p+1$ where the new vertices of the cliques are all different (see Figure~\ref{fig-0}(a)).

Erd\H{o}s, F\"{u}redi, Gould and Gunderson \cite{Inter-tri} determined, for sufficiently large $n$, the extremal number for triangles intersecting in exactly one common vertex. One can think of this graph as blowing up edges of a star to triangles. More generally, $k$ cliques of size $p+1$ intersecting in exactly one common vertex is $S_k^{p+1}$, Chen, Gould, Pfender and Wei \cite{Inter-cli} generalized the main result of \cite{Inter-tri} to $S_k^{p+1}$:

\begin{theorem}\label{inter-clique}\cite{Inter-cli}
For any $p\ge 2$ and $k\ge 1$, and for any $n\ge 16k^3(p+1)^8$, we have
$$\ex(n,S_k^{p+1})=\ex(n,K_{p+1})+
\left\{
  \begin{array}{l l}
    k^2-k & \quad \text{if $k$ is odd,}\\
    k^2-\frac{3}{2}k & \quad \text{if $k$ is even.}\\
  \end{array} \right.$$
\end{theorem}

Given two vertex-disjoint graphs $H$ and $G$, denote by $H\bigotimes G$ the graph obtained by joining each vertex of $H$ to each vertex of $G$. Let $H(n,p,s)$ be $K_{s-1}\bigotimes T_{n-s+1,p}$ (see Figure~\ref{fig-0}(b)) and $H'(n,p,s)$ be any of the graphs obtained by putting one extra edge in any class of $T_{n-s+1,p}$ in $H(n,p,s)$.

Recently, Glebov \cite{Cli-path} determined, for sufficiently large $n$, the extremal number and the extremal graphs for the blow-up of paths. More history of this topic are given in Section~\ref{prel}.

\begin{theorem}\label{cli-path}\cite{Cli-path}
For any $p\ge 2$ and $k\ge 1$, and for any $n>16k^{11}(p+1)^8$, $H(n,p,\lfloor\frac{k-1}{2}\rfloor+1)$ ($H'(n,p,\lfloor\frac{k-1}{2}\rfloor+1)$ resp.) is the unique extremal graph for $P_{k+1}^{p+1}$ when $k$ is odd (even resp.).
\end{theorem}

The main motivation for this note is that \cite{Inter-cli}, \cite{Inter-tri} and \cite{Cli-path} give sporadic results about problems of the same flavor. We unite these extremal problems for blow-ups of graphs and look at the general theory behind these results by investigating the decomposition families of the forbidden graphs. Using the method in \cite{S-sym} (see also \cite{S-I}, \cite{S-H}), we determine the extremal number and found all extremal graphs for all blow-ups of cycles. Somewhat surprisingly, the result for blow-ups of cycles is not much different from blow-ups of paths except for $C_3^3$. Before stating our results, we need a definition. Let $H^*(n)$ be graphs obtained by putting (almost) perfect matchings in both classes in $K_{\lceil n/2\rceil,\lfloor n/2\rfloor}$ (see Figure~\ref{fig-0}(c)).

\begin{theorem}\label{cli-cycle}
For any $p\ge 2$ and $k\ge 3$, when $n$ is sufficiently large, we have the following results:

(i) For any $C_k^{p+1}\neq C_3^3$, $H(n,p,\lfloor\frac{k-1}{2}\rfloor+1)$ ($H'(n,p,\lfloor\frac{k-1}{2}\rfloor+1)$ resp.) is the unique extremal graph for $C_{k}^{p+1}$ when $k$ is odd (even resp.).

(ii) For $C_3^3$, if $4|n$, $H^*(n)$ is the unique extremal graph; otherwise both $H^*(n)$ and $H(n,2,2)$ are extremal graphs.
\end{theorem}

\begin{figure}[ht]
\begin{center}
\includegraphics{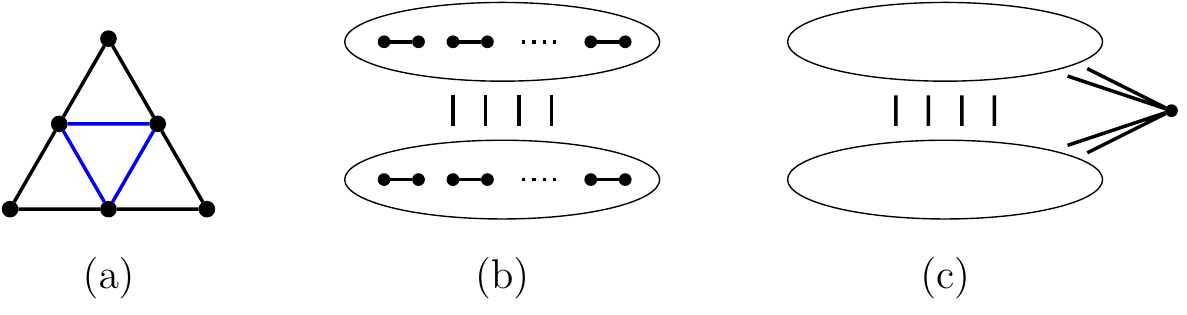}
\end{center}
\caption{(a) $C_3^3$; (b) $H^*(n)$; (c) $H(n,2,2)$.}
\label{C33}
\end{figure}

\begin{rmk}
Both $H'(n,p,s)$ and $H^*(n)$ might contain non-isomorphic graphs. However graphs in the same family are similar in the sense that they have the same number of edges. Since their difference does not matter in this note, we will treat each of them as a ``unique'' graph instead of families of graphs.
\end{rmk}

In addition, for a large class of trees, we determined the extremal number for their blow-up graphs and found their unique extremal graph.

\begin{theorem}\label{cli-trees}
Given a tree $T$, denote by $A$ and $B$ its two color classes with $|A|\le |B|$. For any $p\ge 3$, when $n$ is sufficiently large, we have that

(i) if $T$ has a leaf in $A$ and $\alpha(T)=|B|$, then $H(n,p,|A|)$ is the unique extremal graph for $T^{p+1}$.

(ii) if the minimum degree in $A$ is 2, then $H'(n,p,|A|)$ is the unique extremal graph for $T^{p+1}$.
\end{theorem}

\begin{rmk}
Trees considered in Theorem~\ref{cli-trees} (i) include even paths and those in (ii) include odd path. Hence it implies Theorem~\ref{cli-path} when $p\ge 3$. For $p=2$, the technique in the proof of Theorem~\ref{blow-up cycle} (see Appendix) works for blow-ups of paths. It is not difficult to see that in a proper subdivision of any star, its smaller color class either has minimum degree 2 or has a leaf and its independence number equals to the size of the larger color class, thus Theorem~\ref{cli-trees} can be applied to blow-ups of a proper subdivision of stars, which is an extention of Theorem~\ref{inter-clique}.
\end{rmk}

The rest of this paper is organized as follows: in Section~\ref{prel} we provide more motivation and the key lemma. Section~\ref{blow-up cycle} gives a proof for Theorem~\ref{cli-cycle} when $p\ge 3$ and Section~\ref{blow-up tree} is devoted to proof of Theorem~\ref{cli-trees} (i). The proof for Theorem~\ref{cli-trees} (ii) is similar, we include a sketch of its proof in Appendix together with a proof for Theorem~\ref{cli-cycle} when $p=2$.

We finish this section with a few more definitions that will be used later. We denote the \emph{degree} of a vertex $v$ by $d(v)$ and write $N(v)$ for the set of its neighbors and for $S\subseteq V(G)$, let $N(S)$ be the set of vertices that have some neighbors in $S$. Denote by $K_t^-$ the graph obtained from deleting an edge from a complete graph on $t$ vertices.
A \emph{dominating vertex} in $G$ is a vertex that is adjacent to all other vertices in $G$.
A \emph{linear forest} is a forest whose connected components are paths.
Two disjoint vertex sets $U$ and $W$ are \emph{completely joined} in $G$ if $uw\in E(G)$ for all $u\in U, w\in W$. Write $kH$ for the vertex disjoint union of $k$ copies of $H$. For two vertex-disjoint graphs $H$ and $G$, denote by $H\cup G$ the disjoint union of $H$ and $G$.


\section{Motivation and History}\label{prel}

Given a graph $H$, a \emph{vertex split} on some vertex $v\in V(H)$ is defined as follows: replace $v$ by an independent set of size $d(v)$ in which each vertex is adjacent to exactly one distinct vertex in $N_H(v)$. Given a vertex subset $U\subseteq V(H)$, a vertex split on $U$ means applying vertex split on the vertices in $U$ one by one. It is not difficult to see that the order of vertices we apply vertex split does not matter. Denote by $\mathcal{H}(H)$ the family of graphs that can be obtained from $H$ by applying vertex split on some $U\subseteq V(H)$. Note that $U$ could be empty, therefore $H\in \mathcal{H}(H)$. For example, $\mathcal{H}(P_{k+1})$ is the family of all linear forests with $k$ edges and $\mathcal{H}(C_k)$ consists of $C_k$ and all linear forests with $k$ edges. Given a family $\mathcal{L}$, define $p=p(\mathcal{L})=\displaystyle\min_{L\in \mathcal{L}}\chi(L)-1$.

\begin{defn}
Denote by $I_v$ the $v$-vertex graph with no edges. Given a family $\mathcal{L}$, let $\mathcal{M}:=\mathcal{M}(\mathcal{L})$ be the family of minimal graphs $M$ for which there exist an $L\in\mathcal{L}$ and a $t=t(L)$ such that $L\subseteq M'\bigotimes K_{p-1}(t,\ldots,t)$, where $M'=M\cup I_t$. We call $\mathcal{M}$ the \emph{decomposition family} of $\mathcal{L}$.
\end{defn}

Thus, a graph $M$ is in $\mathcal{M}$ if the graph obtained from putting $M$ into a class of a large $T_{n,p}$ contains some $L\in\mathcal{L}$. If $L\in\mathcal{L}$ with minimum chromatic number $p+1$, then $L\subseteq K_{p+1}(t,\ldots,t)$ for some $t\ge 1$, therefore the decomposition family $\mathcal{M}$ always contains some bipartite graphs.

\begin{expl}
Denote by $O_6$ the edge-graph of the octahedron, namely $O_6=K_{2,2,2}$. Since $O_6=C_4\bigotimes I_2$, we have that $\mathcal{M}(O_6)=\{C_4\}$. For any $\ell>1$, $\mathcal{M}(C_{2\ell+1})=\{P_2\}$. For blow-ups of stars, $\mathcal{M}(S_k^{p+1})=\{S_k, M_k\}$ for $p\ge 2$ and $k\ge 1$; for blow-ups of paths $\mathcal{M}(P_{k+1}^{p+1})=\{$all linear forests with $k$ edges$\}=\mathcal{H}(P_{k+1})$ and for cycles $\mathcal{M}(C_k^{p+1})=\{C_k, $ all linear forests with $k$ edges$\}=\mathcal{H}(C_k)$ for $p\ge 3$ and $k\ge 1$.
\end{expl}

For a family of forbidden graph $\mathcal{L}$ with decomposition family $\mathcal{M}$, we have
\begin{equation}\label{eq:bounds}
e(T_{n,p})+ \ex(\frac{n}{p},\mathcal{M})\le \ex(n,\mathcal{L})\le e(T_{n,p})+(1+o(1))p\cdot \ex(\frac{n}{p},\mathcal{M}),
\end{equation}
where the lower bound is obtained from putting an $\mathcal{M}$-extremal graph in one of the classes in $T_{n,p}$ (see \cite{bib-bollobas-EGT}).

The Erd\H{o}s-Stone-Simonovits Theorem determines asymptotically the extremal functions of non-bipartite graphs, while the decomposition family governs the finer error terms as shown in~\eqref{eq:bounds}, thus it helps to give sharper bounds on the extremal number.

There are examples where the upper bound in~\eqref{eq:bounds} holds. Let $Q(r,p)$ be the graph consisting of a dominating vertex and a $p$-class Tur\'an graph on $rp$ vertices, namely $Q(r,p)=T_{rp,p}\bigotimes I_1$. Notice that $\mathcal{M}(Q(r,p))=\{S_r\}$, thus an $\mathcal{M}(Q(r,p))$-free graph has maximum degree $r-1$. Simonovits \cite{Simonovits1966} showed that a $Q(r,p)$-extremal graph can be obtained from putting (almost) $(r-1)$-regular triangle-free graphs into each class of a $T_{n,p}$.

We recall the Octahedron Theorem by Erd\H{o}s and Simonovits \cite{octahedron}, which gives an example where neither the upper bound nor the lower bound in~\eqref{eq:bounds} is true.

\begin{theorem}\label{oct}
For sufficiently large $n$, every $O_6$-extremal graph $S_n$ can be obtained as $S_n=U_m\bigotimes Z_{n-m}$, for some $m$-vertex $C_4$-extremal graph $U_m$ and some $(n-m)$-vertex $P_3$-extremal graph $Z_{n-m}$, where $m=n/2+o(n)$.
\end{theorem}

A graph $L$ is \emph{weakly edge-color-critical}, or shortly \emph{weakly-critical}, if there is an edge $e\in E(L)$ for which $\chi(L-e)<\chi(L)$. Simonovits \cite{Simonovits1966} proved that the Tur\'an graph is the unique extremal graph for weakly-critical graphs when $n$ is sufficiently large. In the same paper, he also proved when the forbidden graph is $L=sH$, where $H$ is weakly-critical, and $\chi(H)=p+1\ge 3$, then for sufficiently large $n$, the unique $L$-extremal graph is $H(n,p,s)$. Later in \cite{S-sym}, he further generalized this result to the following theorem.

\begin{theorem}\label{Hnps-thm}
Let $\mathcal{L}$ be the family of forbidden graphs and $p=p(\mathcal{L})=\displaystyle\min_{L\in \mathcal{L}}\chi(L)-1$. If by omitting any $s-1$ vertices of any $L\in \mathcal{L}$ we obtain a graph with chromatic number at least $p+1$, but by omitting $s$ suitable edges of some $L\in \mathcal{L}$ we get a $p$-colorable graph, then $H(n,p,s)$ is the unique extremal graph for $n$ sufficiently large.
\end{theorem}

Simonovits~\cite{when-Hnps} asked the following question.
\begin{ques}\label{question}
Characterize graphs whose unique extremal graph is of the form $H(n,p,s)$.
\end{ques}

We make a step towards answering Question~\ref{question}: notice that the blow-ups of cycles and trees does not satisfy the hypothesis in Theorem~\ref{Hnps-thm}, hence Theorem~\ref{cli-cycle} and~\ref{cli-trees} provides an additional family of forbidden graphs whose unique extremal graph is $H(n,p,s)$ for suitable $p$ and $s$.

On the other hand, the results in \cite{Inter-cli}, \cite{Inter-tri}, \cite{Cli-path}, \cite{Simonovits1966} and \cite{S-sym} show that for blow-ups of stars, paths and many other families of graphs, the lower bound construction is optimal. Our results show that for blow-ups of cycles and a large class of trees, this is also the case. It would be interesting to describe all decomposition families $\mathcal{M}$ where the lower bound in~\eqref{eq:bounds} is sharp. Here we make the first attempt towards this direction.

The following lemma shows that the decomposition family of blow-ups of some graphs (in particular bipartite graphs) is actually the family obtained from splitting its vertices.

\begin{lemma}\label{key-lem}
Given $p\ge 3$ and any graph $H$ with $\chi(H)\le p-1$, $\mathcal{M}(H^{p+1})=\mathcal{H}(H)$. In particular, a matching of size $e(H)$ is in $\mathcal{M}(H^{p+1})$.
\end{lemma}

\noindent\textbf{Proof.}
Note that $\chi(H^{p+1})=p+1$, by the definition of decomposition family, any graph $M$ in $\mathcal{M}(H^{p+1})$ is a minimal graph under the condition that there exists a copy of $M$ in $H$, such that removing the vertex set of $M$ together with a suitable independent set $I_v$ results in a $(p-1)$-colorable graph. Since $I_v$ is an independent set, it can have at most one vertex from each $(p+1)$-clique in $H^{p+1}$. Recall the removal of $M\cup I_v$ decreases the chromatic number by at least two, this implies $M$ should have at least one vertex from each $(p+1)$-clique in $H^{p+1}$. Since a vertex from $M$ and a vertex from $I_v$ in the same $(p+1)$-clique would be adjacent, which contradicts the minimality of $M$, thus $I_v$ is empty and $M$ includes exactly two vertices from each $(p+1)$-clique in $H^{p+1}$. Locally, for each $v\in H$, $N_H(v)\cup \{v\}$ spans a $S_{d_H(v)}^{p+1}$ in $H^{p+1}$.
If the center of that $S_{d_H(v)}^{p+1}$ is in $M$, then $M$ also contains exactly one other vertex from each of the $(p+1)$-cliques in this $S_{d_H(v)}^{p+1}$. This implies $v$ is not a split vertex in $M$. Otherwise, $M$ contains two vertices from each $(p+1)$-cliques of $S_{d_H(v)}^{p+1}$, which implies that $v$ was split into $d_H(v)$ leaves in $M$. Hence $\mathcal{M}(H^{p+1})=\mathcal{H}(H)$. In particular, by splitting all vertices of $H$, we obtain a matching of size $e(H)$.
\qed

\medskip

The following definition was introduced in \cite{S-sym}.

\begin{defn}\label{def:sym}
Denote by $\mathbb{D}(n,p,r)$ the family of $n$-vertex graphs $G_n$ satisfying the following symmetry condition:

(i) It is possible to omit at most $r$ vertices of $G_n$ so that the remaining graph $G'$ is a product of graphs of almost equal order:
$G'=\displaystyle\prod_{i\le p}G^i$, where $|V(G^i)|=n_i$ and $|n_i-\frac{n}{p}|\le r$, for every $i\le p$.

(ii) For every $i\le p$, there exist connected graphs $H_i$ such that $G_i=k_iH_i$, where $k_i=\frac{n_i}{|V(H_i)|}$ and any two copies $H_i^j,H_i^{\ell}$ in $G_i$ ($1\le j< \ell\le k_i$), are \emph{symmetric subgraphs} of $G_n$: there exists an isomorphism $\omega:H_i^j\rightarrow H_i^{\ell}$ such that for every $x\in H_i^j$, $u\in G_n-H_i^j-H_i^{\ell}$, $xu\in E(G_n)$ if and only if $\omega(x)u\in E(G_n)$.
\end{defn}

\begin{figure}[ht]
\begin{center}
\includegraphics{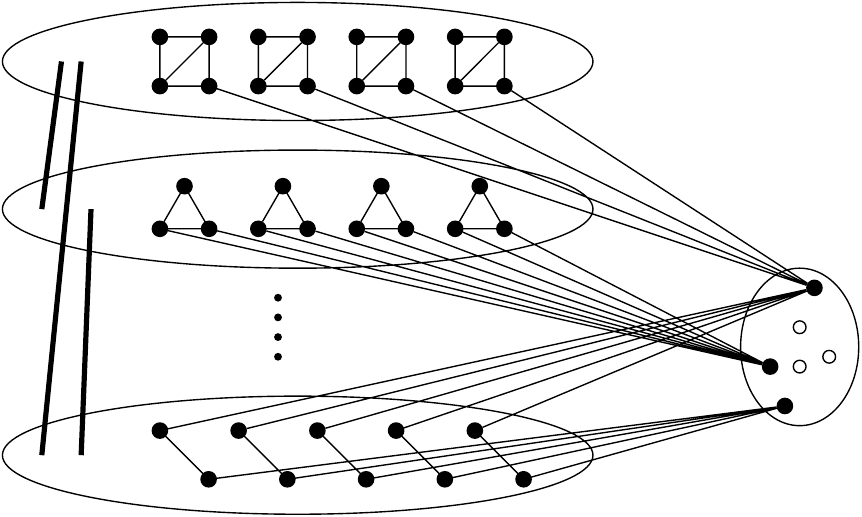}
\end{center}
\caption{Symmetric subgraphs: blocks here are $K_4^-$, $K_3$ and $P_2$ respectively.}
\label{symm}
\end{figure}

The graphs $H_i$ will be called the \emph{blocks}, the vertices in $G_n-G'$ will be called \emph{exceptional} vertices (see Figure~\ref{symm}).

We will need the following two results of Simonovits (\cite{S-sym, S-H}).

\begin{theorem}\label{th:S-sym}\cite{S-sym}
Assume that a finite family $\mathcal{L}$ of forbidden graphs with $p(\mathcal{L})=p$ is given. If for some $L\in \mathcal{L}$ and $\ell:=|V(L)|$,
\begin{equation}\label{eq:path}
L\subseteq P_{\ell}\bigotimes K_{p-1}(\ell,\ell,\ldots,\ell),
\end{equation}
then there exist $r=r(L)$ and $n_0=n_0(r)$ such that $\mathbb{D}(n,p,r)$ contains an $\mathcal{L}$-extremal graph for every $n>n_0$. Furthermore, if this is the only extremal graph in $\mathbb{D}(n,p,r)$, then it is the unique extremal graph for every sufficiently large $n$.
\end{theorem}

\begin{theorem}\label{S-H}\cite{S-H}
Assume that a finite family $\mathcal{L}$ of forbidden graphs with $p(\mathcal{L})=p$ is given. If for some $L\in \mathcal{L}$ and $\ell:=|V(L)|$,
\begin{equation}\label{eq:matching}
L\subseteq \ell P_2\bigotimes K_{p-1}(2\ell,2\ell,\ldots,2\ell),
\end{equation}
then there exist $r=r(L)$ and $n_0=n_0(r)$ such that $\mathbb{D}(n,p,r)$ contains an $\mathcal{L}$-extremal graph for every $n>n_0$. Furthermore for any $\mathcal{L}$-extremal graph $G\in\mathbb{D}(n,p,r)$, we have that

(i) all blocks of $G$ will consist of isolated vertices: the product graph $G'$ will be a Tur\'{a}n graph $T_{n',p}$.

(ii) each exceptional vertex in $G-G'$ is joined either to all the vertices of $G'$ or to all the vertices of $p-1$ classes of $G'$ and to no vertex of the remaining class.
\end{theorem}

The key idea of the proof of Theorem~\ref{cli-cycle} is using Theorem~\ref{S-H} to get a good vertex partition of an extremal graph of $C_k^{p+1}$. Then show that in this partition, there are $t$ exceptional vertices, where $t=\lfloor\frac{k-1}{2}\rfloor$, and the remaining of the graph is a Tur\'an graph (with one extra edge if $k$ is even). This together with $G$ being extremal imply $G\simeq H(n,p,t+1)$ ($H'(n,p,t+1)$ if $k$ is even).

\begin{rmk}
Assume first $k$ is odd. Notice that $H(n,p,t+1)$ is $C_{k}^{p+1}$-free. Indeed the set of $t$ dominating vertices in $H(n,p,t+1)$ together with one class of $T_{n-t,p}$ is $\mathcal{H}(C_{k})$-free. By Lemma~\ref{key-lem}, $\mathcal{M}(C_k^{p+1})=\mathcal{H}(C_{k})$ when $p\ge 3$, and observe that when $p=2$, $\mathcal{M}(C_k^{3})\subseteq \mathcal{M}(C_k^{p+1})$. Thus $H(n,p,t+1)$ is the product of $p-1$ empty graphs and one $\mathcal{M}(C_{k}^{p+1})$-free graph, which is the lower bound construction in~\eqref{eq:bounds}. Hence $H(n,p,t+1)$ is $C_{k}^{p+1}$-free. Similarly when $k$ is even, $H'(n,p,t+1)$ is $C_{k}^{p+1}$-free.
\end{rmk}


\section{Proof of Theorem~\ref{cli-cycle}}\label{blow-up cycle}

\begin{figure}[ht]
\begin{center}
\includegraphics{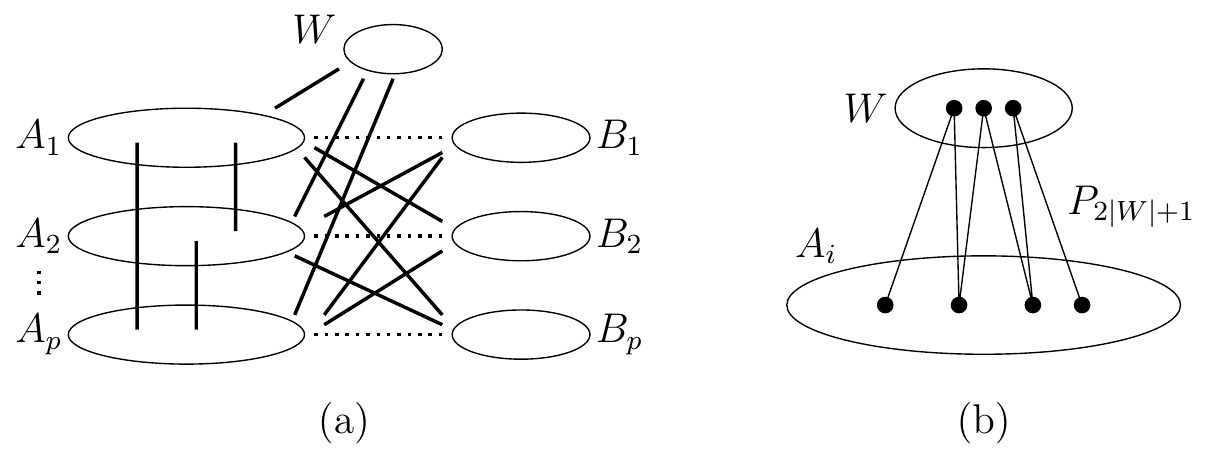}
\end{center}
\caption{(a): a partition of $G$; (b): If $|W|>t$, then $2|W|+1\ge k+1$, $P_{k+1}\subseteq G[W\cup A_i]$.}
\label{longpath}
\end{figure}

\noindent\textbf{Proof.}
We prove Theorem~\ref{cli-cycle} for $p\ge 3$. A crucial observation is that \eqref{eq:path} is equivalent to $\mathcal{M}(C_k^{p+1})$ containing some linear forest of size at most $\ell$ and (3) is equivalent to $\mathcal{M}(C_k^{p+1})$ containing a matching of size at most $\ell$. Since $P_{k+1}\in\mathcal{M}(C_k^{p+1})$, \eqref{eq:path} is satisfied with $L=C_k^{p+1}, \ell=k+1$. Thus by Theorem~\ref{th:S-sym}, there exists an extremal graph $G$ in $\mathbb{D}(n,p,r)$ of $C_{k}^{p+1}$ for some $r$. It suffices to prove that $G\simeq H(n,p,t+1)$ ($H'(n,p,t+1)$ resp.) when $k$ is odd (even resp.). Then Theorem~\ref{th:S-sym} implies it would be the unique extremal graph. Since $M_k\in\mathcal{M}(C_k^{p+1})$, \eqref{eq:matching} is also satisfied with $L=C_k^{p+1}, \ell=k$. We can apply Theorem~\ref{S-H} to get a vertex partition of $G$. Let $A_1,\ldots,A_p$ be the $p$ classes in $T_{n',p}$. Let $W$ be the set of vertices in $G-T_{n',p}$ that are joined to all vertices in $T_{n',p}$ and let $B_i$ be the set of vertices in $G-T_{n',p}-W$ that are joined to all the vertices in $T_{n',p}$ but $A_i$ (see Figure~\ref{longpath}(a)). Define $C_i=A_i\cup B_i$, for all $i$. Note that in $G$ all the cross-edges between $A_i$ and $C_j$ with $i\neq j$ are present, there might be some missing edges between some $B_i$ and $B_j$. Let $D_i\subseteq C_i$ consist of vertices with no neighbor in $W$. Recall that $\mathcal{M}(C_k^{p+1})=\{C_k, $ all linear forests with $k$ edges$\}$. We will frequently use the following fact.

\begin{claim}\label{M}
For any $i\le p$, $G[W\cup C_i]$ is $\mathcal{M}(C_k^{p+1})$-free.
\end{claim}
\noindent\textbf{Proof.}
Notice that $W\cup C_i$ is completely joined to $\displaystyle\bigcup_{j\neq i}A_j$. If for some $M\in\mathcal{M}(C_k^{p+1})$, $M\subseteq G[W\cup C_i]$, then $M'\bigotimes K_{p-1}(n'/p,\ldots,n'/p)\subseteq G$. This, by the definition of decomposition family, implies $C_k^{p+1}\subseteq G$, a contradiction.
\qed
\medskip

\begin{claim}\label{c1}
$|W|=t=\lfloor\frac{k-1}{2}\rfloor$.
\end{claim}
\noindent\textbf{Proof.}
First of all, $|W|\le t$. Indeed, by Claim~\ref{M}, $G[W\cup C_i]$ is $\mathcal{M}(C_k^{p+1})$-free. Since $W$ and $A_i$ are completely joined, if $|W|>t$, then $P_{k+1}\subseteq G[W\cup A_i]$ (see Figure~\ref{longpath}(b)). But $P_{k+1}\in\mathcal{M}(C_k^{p+1})$, a contradiction.

On the other hand, suppose $|W|\le t-1$, then some simple calculation shows
$$e(G)\le e(T_{n,p})+\frac{t-1}{p}n+o(n).$$
However since $C_{k}^{p+1}\not\subseteq H(n,p,t+1)$, we have
$$e(G)\ge e(H(n,p,t+1)) \ge e(T_{n,p})+\frac{t}{p}n+o(n),$$
a contradiction.
\qed
\medskip

\noindent\textbf{Case 1:} $k$ is odd. Then $t=\frac{k-1}{2}$.
We shall show that for each $i\le p$, $G[C_i]$ has no edge. Then by Claim~\ref{c1} and the maximality of $G$,  $G\simeq H(n,p,t+1)$. Indeed, any edge $xy\in C_i$ together with a $P_{2|W|+1}=P_k$ in $G[W\cup A_i]$ avoiding $\{x,y\}$ form a linear forest with $k$ edges, which is in $\mathcal{M}(C_{k})$. This contradicts Claim~\ref{M}.

\medskip

\noindent\textbf{Case 2:} $k$ is even. Then $t=k/2-1$. It suffices to prove following claim. This together with Claim~\ref{c1} and the maximality of $G$ implies $G\simeq H'(n,p,t+1)$.

\begin{claim}\label{c2}
When $k$ is even, we have (i) each $G[C_i]$ has at most one edge; (ii) the edge in $C_i$ has at least one endpoint not in $D_i$; and (iii) there is at most one class $G[C_i]$ having one such edge.
\end{claim}

\noindent\textbf{Proof.}
(i): For contradiction suppose there are two edges in some $G[C_i]$, say $e_1$ and $e_2$. Then one can find a copy of $P_{2|W|+1}=P_{k-1}$ in $G[W\cup A_i]$ avoiding the endpoints of $e_1,e_2$, so we get a linear forest with $k$ edges in $G[W\cup C_i]$, a contradiction.

(ii): Suppose $u_iv_i\in E(G[C_i])$ and $u_i,v_i\in D_i$. Since $k$ is even hence $k\ge 4$, therefore $W$ is not empty: $|W|=t=k/2-1\ge 1$. Denote $G_0$ the graph obtained by deleting the edge $u_iv_i$ and adding all edges between $\{u_i,v_i\}$ and $W$. There are at least two such cross-edges since $W$ is nonempty and $u_i,v_i\in D_i$. Thus $e(G_0)>e(G)$. It remains to show $G_0$ is also $C_{k}^{p+1}$-free, which contradicts to the extremality of $G$. Notice that $u_i$ and $v_i$ is not joined to any vertex in $C_i$ in $G_0$. Hence they have the same neighborhood as vertices in $A_i$. Also since $G$ and $G_0$ only differ at $u_i$ and $v_i$, a copy of $C_k^{p+1}$ in $G_0$ must involve $u_i$ or $v_i$ or both. But then one can obtain a copy of $C_k^{p+1}$ in $G$ by replacing vertices from $\{u_i,v_i\}$ by vertices from $A_i$, a contradiction.

\begin{figure}[ht]
\begin{center}
\includegraphics{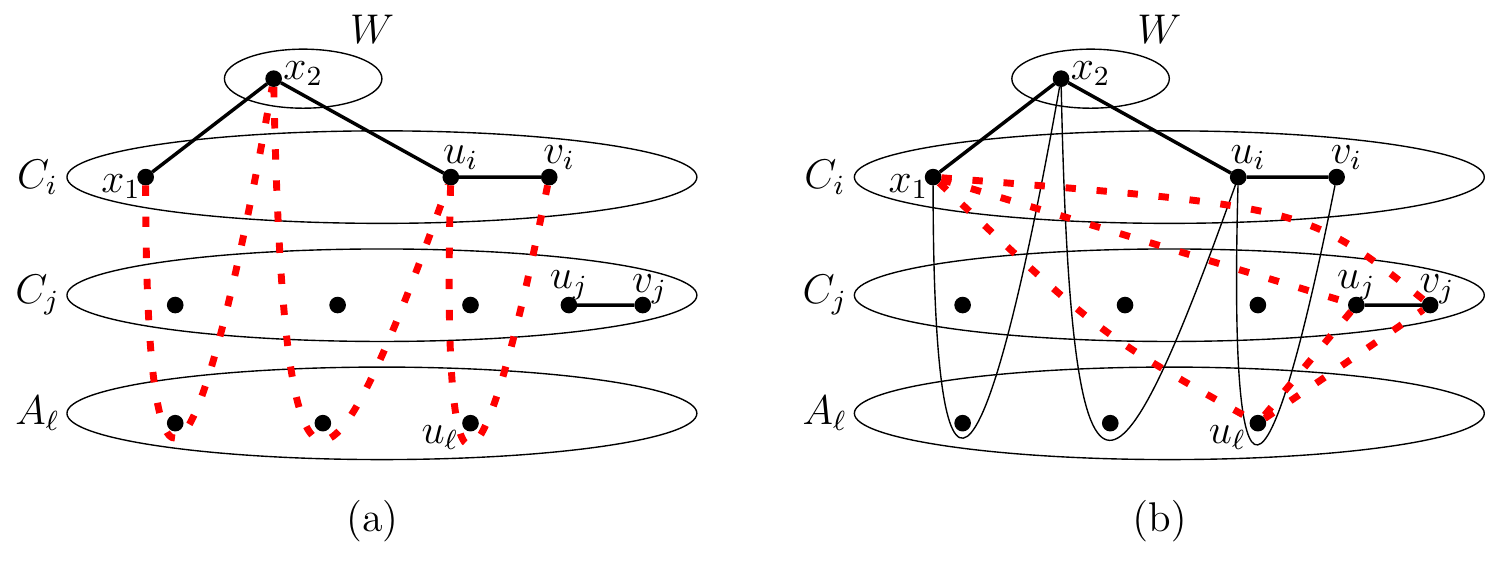}
\end{center}
\caption{When $k=4$: (a) Extend a $P_k$ to a $P_k^{p+1}$; (b) Obtain a $C_{k}^{p+1}$ from $P_{k}^{p+1}$.}
\label{path-cyl}
\end{figure}

(iii): Suppose for some $i\neq j$, $u_iv_i\in E(G[C_i])$ and $u_jv_j\in E(G[C_j])$ with $u_i\not\in D_i$ and $u_j\not\in D_j$. Then we can find a copy of $P_k$ in $G[C_i\cup W]$ which starts at a vertex in $A_i$ and whose last edge is $u_iv_i$. Denote vertices of such a path $x_1,x_2,\ldots,x_{k-2},u_i,v_i$ with $x_1\in A_i$. We can then extend this path to a $P_k^{p+1}$ in $G$ (see Figure~\ref{path-cyl}(a)).


Since $p\ge 3$, there is a third class, say $A_{\ell}$ with $\ell\neq i,j$. The last clique on this copy of $P_k^{p+1}$, namely the one containing $u_i$ and $v_i$, intersect $A_{\ell}$ at exactly one vertex, call it $u_{\ell}$. Then this $P_k^{p+1}$ together with the $(p+1)$-clique consisting of one vertex from each $A_q$, $q\neq i,j,\ell$, and $u_{\ell}, x_1,u_j,v_j$, form a $C_{k}^{p+1}$ (see Figure~\ref{path-cyl}(b)), where $x_1,\ldots,x_{k-2},u_i,u_{\ell},x_1$ is the vertices of $C_k$ that was blown up. This yields a contradiction.
\qed


\section{Proof of Theorem~\ref{cli-trees}}\label{blow-up tree}
In this section, unless otherwise specified, $p\ge 3$ and $T$ is a tree with two color classes (partite sets) $A$ and $B$ such that $|A|\le |B|$. Let $a=|A|$, $b=|B|$. Recall that, by Lemma~\ref{key-lem}, $\mathcal{M}(T^{p+1})=\mathcal{H}(T)$. In particular, $T\in \mathcal{M}(T^{p+1})$.

\begin{lemma}\label{T-free}
If $T$ has a leaf in $A$ and $\alpha(T)=b$, then for any $m\ge 1$, $K_{a-1}\bigotimes \overline{K}_m$ is $\mathcal{H}(T)$-free, hence $\mathcal{M}(T^{p+1})$-free.
\end{lemma}

\begin{figure}[ht]
\begin{center}
\includegraphics{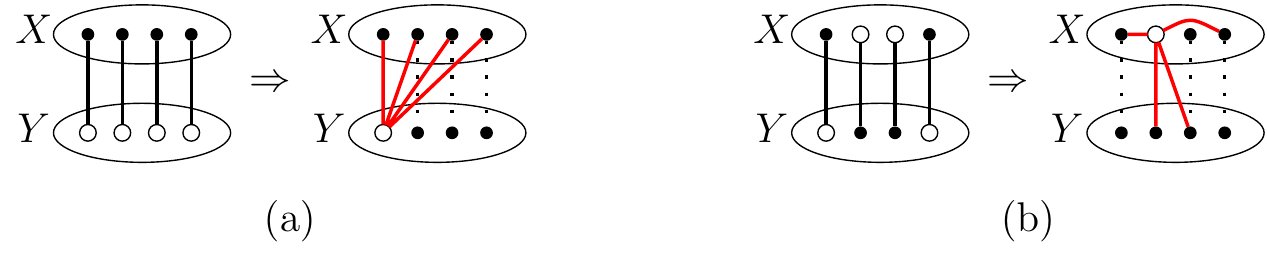}
\end{center}
\caption{Circled vertices are in $L(v)$. (a): $L(v)\subseteq Y$; (b): $L(v)\cap X\neq \emptyset$.}
\label{figTfree}
\end{figure}

\noindent\textbf{Proof.}
For simplicity, let $G=K_{a-1}\bigotimes \overline{K}_m$ with $V(G)=X\cup Y$, where $X$ is the set of vertices in the $(a-1)$-clique and $Y$ is the remaining independent set. We may assume $|Y|=m\ge b+1$, otherwise $|V(G)|<|V(F)|$ for any $F\in \mathcal{H}(T)$. First notice that $T\not\subseteq G$. Indeed, an embedding of $T$ in $G$ has at least $|V(T)|-(a-1)=b+1$ vertices in $Y$, which contradicts to $\alpha(T)=b$. For contradiction, suppose that some forest $F\in \mathcal{H}(T)$ is in $G$. Recall that $F$ is obtained from splitting vertices in some $U\subseteq V(T)$ to sets of leaves in $F$. For any $v$ in $U$, denote by $L(v)\subseteq V(G)$ the set of leaves in $F$ corresponding to $v$. We shall get a copy of $T$ in $G$ from a copy of $F$ by applying the following operation to every $v\in U$ to undo the vertex split: pick any $v\in U$, look at the corresponding $L(v)$ in $F$. First discard edges in $F$ adjacent to $L(v)$. Then if $L(v)\subseteq Y$, add in $F$ edges from $N_G(L(v))$ to a vertex in $Y$; otherwise add in $F$ edges from $N(L(v))$ to a vertex in $L(v)\cap X$ (see Figure~\ref{figTfree}).
\qed
\medskip

From Lemma~\ref{T-free}, we immediately get the following.

\begin{prop}
$H(n,p,a)$ is $T^{p+1}$-free.
\end{prop}

\noindent\textbf{Proof.}
Since $H(n,p,a)$ is the product of $p-1$ independent sets and one $\mathcal{M}(T^{p+1})$-free graph (obtained from combining the set of $a-1$ dominating vertices with the remaining independent set). This is the lower bound construction in~\eqref{eq:bounds}. Hence $H(n,p,a)$ is $T^{p+1}$-free.
\qed
\medskip

\noindent\textbf{Proof of Theorem~\ref{cli-trees} (i).} Since a matching of size $e(T)$ is in $\mathcal{M}(T^{p+1})$, we can proceed as in the proof of Theorem~\ref{cli-cycle} and define $W,A_i,B_i,C_i$ in the same way. If $|W|\ge a$, then $T\subseteq G[W\cup C_i]$ for any $i$. However $T\in \mathcal{M}(T^{p+1})$, this implies $T^{p+1}\subseteq G$, a contradiction. If $|W|\le a-2$, then
$$e(G)=e(T_{n,p})+\frac{a-2}{p}n+o(n)<e(T_{n,p})+\frac{a-1}{p}n+o(n)=e(H(n,p,a)),$$
a contradiction. Thus $|W|=a-1$. Also we may assume that $W$ is non-empty. Indeed, if $W=\emptyset$, then $a=1$. Since $T$ has a leaf in $A$, it implies $T$ is $P_2$, then $T^{p+1}$ is $K_{p+1}$ and its unique extremal graph is $H(n,p,1)=T_{n,p}$.

It remains to show that $e(G[C_i])=0$ for all $1\le i\le p$. Let $u$ be a leaf of $T$ in $A$, and $v$ be its neighbor in $B$. Let $T'=T-u$ and let $F$ and $F'$ be the forests obtained from splitting $v$ in $T$ and $T'$ respectively. Notice that $F=F'\cup K_2$. Indeed, $v$ has one more neighbor (leaf $u$) in $T$, which becomes a $K_2$ after splitting. Suppose there is an edge $xy$ in $G[C_i]$ for some $i$. Since $u\in A$ and $|W|=a-1$, $T'$ has an embedding in $G[W\cup C_i- \{x,y\}]$ with $v$ in $C_i$. Splitting $v$ in this copy of $T'$ (not using $x$ or $y$), we get a copy of $F'$ in $W\cup C_i-\{x,y\}$. Note that for any $i\le p$, $G[W\cup C_i]$ is $\mathcal{M}(T^{p+1})$-free, since $W\cup C_i$ is completely joined to $A_j$, $\forall j\neq i$. Thus edge $xy$ together with this $F'$ yields a copy of $F$ in $G[W\cup C_i]$. However $F\in \mathcal{H}(T)=\mathcal{M}(T^{p+1})$, a contradiction.
\qed
\section*{Acknowledgments}
The author would like to thank J\'ozsef Balogh for encouragement and valuable comments and remarks. He would also like to thank Roman Glebov for helpful discussions.




\appendix
\section*{Appendix}

\begin{figure}[ht]
\begin{center}
\includegraphics{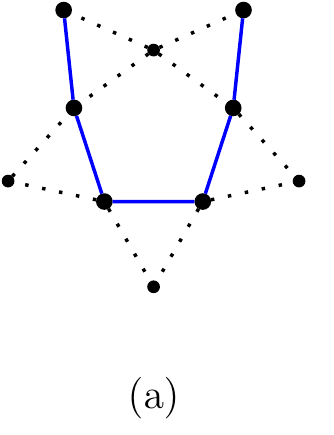}
\hspace{2cm}
\includegraphics{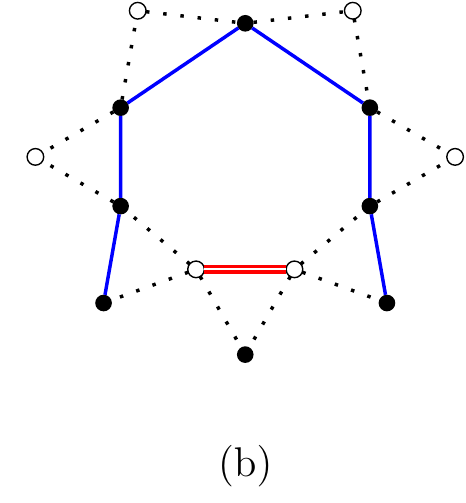}
\end{center}
\caption{(a): $P_{k+1}\in\mathcal{M}(C_k^3)$ when $k=5$; (b): $C_k^3\subseteq (P_k\cup I_k)\bigotimes(P_2\cup I_k)$ when $k=7$.}
\label{fig-1}
\end{figure}

\noindent\textbf{Proof of Theorem~\ref{cli-cycle} for $p=2$.}
Given $C_k^3$, notice that $P_{k+1}\in \mathcal{M}(C_k^3)$ (See Figure~\ref{fig-1}(a)). Thus by Theorem~\ref{th:S-sym}, there is an extremal graph $G\in \mathbb{D}(n,2,r)$ for $C_k^3$ and removing a few exceptional vertices from $G$ results in $G'=A_1\bigotimes A_2$, where $A_1$ and $A_2$ are disjoint unions of symmetric subgraphs (blocks) $H_1$ and $H_2$ respectively.

\begin{claim}\label{p=2}
For any $k\ge 4$, each block $H_i$, $i=1,2$, is an isolated vertex. For $C_3^3$, $H_i$ can be an isolated vertex or $P_2$, and $H^*(n)$ is the unique extremal graph when $4|n$, otherwise $H(n,2,2)$ is also an extremal graph.
\end{claim}

We first show how Claim~\ref{p=2} implies Theorem~\ref{cli-cycle} for $p=2$. When $k\ge 4$, since $H_1$ and $H_2$ are symmetric with respect to $G$, Claim~\ref{p=2} implies the exceptional vertices are adjacent either to all the vertices in $G'$ or to all the vertices in one class of $G'$ and none of the other class. Similarly let $W$ be the set of vertices adjacent to all vertices of $G'$, and $B_1$, $B_2$ be the sets of vertices joining only vertices in $A_2$ and $A_1$ respectively. Let $C_i=A_i\cup B_i$, for $i=1,2$. Claim~\ref{c1} is still true, namely $|W|=t$.

When $k$ is odd, $t=\frac{k-1}{2}$. It suffices to show $e(C_i)=0$ for $i=1,2$. Indeed, the 2-coloring in Figure~\ref{fig-1}(b) shows that $C_k^3\subseteq(P_k\cup I_k)\bigotimes(P_2\cup I_k)$. Suppose $C_i$ has an edge, since $P_{2|W|+1}=P_k\subseteq W\cup A_{3-i}$, we have $C_k^3\subseteq(W\cup A_{3-i})\bigotimes C_i\subseteq G$, a contradiction.

\begin{figure}[ht]
\begin{center}
\includegraphics{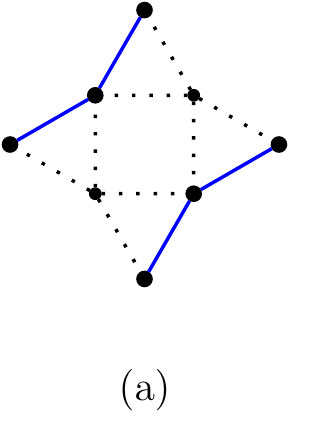}
\hspace{2cm}
\includegraphics{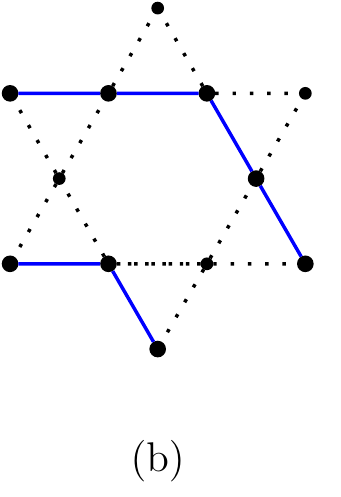}
\end{center}
\caption{(a): $(P_3\cup P_{k-1})\in \mathcal{M}(C_k^3)$ when $k=4$; (b): $(P_5\cup P_{k-3})\in \mathcal{M}(C_k^3)$ when $k=6$.}
\label{fig-2}
\end{figure}

When $k$ is even, $t=k/2-1$. It suffices to show only one class $C_i$ has at most one edge. First notice that each $C_i$ is $P_3$-free, since $(P_3\cup P_{k-1})\in \mathcal{M}(C_k^3)$ (see Figure~\ref{fig-2}(a)) and $P_{2|W|+1}=P_{k-1}\subseteq W\cup A_i$. Suppose some $C_i$ has two isolated edges, say $x_1y_1$ and $x_2y_2$. Then each edge $x_iy_i$ has at least one endpoint adjacent to some vertices in $W$, w.l.o.g. let them be $x_1,x_2$, since otherwise deleting $x_iy_i$ and adding all edges between $W$ and $\{x_i,y_i\}$ (at least two such edges) results in a $C_k^3$-free graph with more edges than $G$, a contradiction.
If $x_1,x_2$ are adjacent to the same vertex $u\in W$, then
$y_1,x_1,u,x_2,y_2$ form a $P_5$, and $P_{k-3}\in (W-u)\cup( A_i-\{x_1,x_2,y_1,y_2\})$. Thus $(P_5\cup P_{k-3})\subseteq (W\cup C_i)$, a contradiction since $(P_5\cup P_{k-3})\in \mathcal{M}(C_k^3)$ (see Figure~\ref{fig-2}(b)). If $x_1,x_2$ are adjacent to different vertices in $W$, then a copy of $P_{k+1}$ can be obtained in $W\cup C_i$ by prolonging a $P_{k-1}$ using edges $x_1y_1$ and $x_2y_2$. Then we get a contradiction since $P_{k+1}\in \mathcal{M}(C_k^3)$. Thus each $C_i$ has at most one edge. Now suppose both $C_i$, $i=1,2$, contain an edge $u_iv_i$ with $u_i$ adjacent to some vertices in $W$. Then similarly we can get a copy of $P_k$ in $G[C_1\cup W]$, starting at a vertex in $A_1$ and ending with edge $u_1v_1$. Let the vertices on this path be $x_1,x_2,\ldots,x_{k-2},u_1,v_1$ with $x_1\in A_1$. We can expand the path $x_1,\ldots,x_{k-2},u_1$ to a copy of $P_{k-1}^3$. Since $x_1\in A_1$, $x_1$ is adjacent to all vertices in $C_2$. In particular, $x_1$ is adjacent to both $u_2$ and $v_2$. Thus $u_2,v_2$ together with that copy of $P_{k-1}^3$ form a $P_{k}^3$. Note that if there are at least three edges between $\{u_1,v_1\}$ and $\{u_2,v_2\}$, then it would complete $P_{k}^3$ to a copy of $C_{k}^3$. Thus there are at most two edges, then delete $u_2v_2$ and add the missing edges between $\{u_1,v_1\}$ and $\{u_2,v_2\}$. The resulting graph is still $C_{k}^3$-free but with more edges than $G$, a contradiction.
\qed
\medskip

\noindent\textbf{Proof of Claim~\ref{p=2}.} We distinguish two cases depending on the parity of $k$.

\begin{figure}[ht]
\begin{center}
\includegraphics{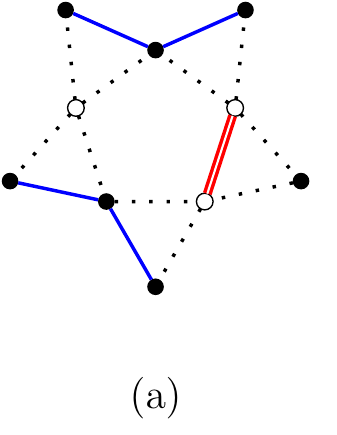}
\hspace{2cm}
\includegraphics{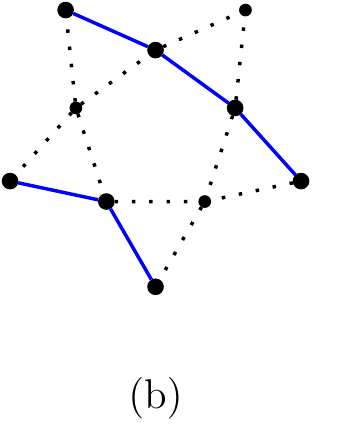}
\end{center}
\caption{(a): $C_k^3\subseteq (\frac{k-1}{2}P_3\cup I_k)\bigotimes (P_2\cup I_k)$ when $k=5$; (b): $(\frac{k-3}{2}P_3\cup P_4)\in\mathcal{M}(C_k^3)$ when $k=5$.}
\label{fig-3}
\end{figure}

\noindent Case 1: $k$ is odd. First we show that $H_i$, $i=1,2$, is $P_3$-free. Suppose to the contrary that $P_3\subseteq H_1$, then $H_2$ has to be an isolated vertex. Since otherwise $\frac{k-1}{2}P_3\cup I_k\subseteq H_1$ and $P_2\cup I_k\subseteq H_2$. This yields a contradiction since $C_k^3\subseteq (\frac{k-1}{2}P_3\cup I_k)\bigotimes (P_2\cup I_k)$ (see Figure~\ref{fig-3}(a)). Furthermore, since $(\frac{k-3}{2}P_3\cup P_4)\in\mathcal{M}(C_k^3)$ (see Figure~\ref{fig-3}(b)), $W$ is empty and $H_1$ is $P_4$-free, otherwise $(\frac{k-3}{2}P_3\cup P_4)\subseteq W\cup C_1$, a contradiction.
If $k=3$, then $t=1$, $K_3\in\mathcal{M}(C_3^3)$, and $H_1$ is $\{K_3,P_4\}$-free, which implies $H_1$ is a star of constant order, say $r$. Then $$e(G)\sim \frac{n^2}{4}+\frac{r-1}{r}\frac{n}{2}< \frac{n^2}{4}+\frac{nt}{p}=\frac{n^2}{4}+\frac{n}{2}\sim e(H(n,2,2)),$$ a contradiction. If $k\ge 5$, then $t\ge 2$. By Erdos-Gallai, since $H_1$ is $P_4$-free, the size of $G$ is maximized when $H_1=K_3$, hence $$e(G)\sim \frac{n^2}{4}+\frac{n}{2}<\frac{n^2}{4}+n\le \frac{n^2}{4}+\frac{nt}{p} \sim e(H(n,2,t+1)),$$ a contradiction. Thus we may assume $H_i$, $i=1,2,$ is $P_3$-free.

If $k=3$ and one of $H_i$ is not an isolated vertex, then $W$ is empty, since otherwise $K_3\subseteq W\cup H_i$. This implies for $i=1,2$, either $H_i=P_1$ and $|W|=1$, or $H_i=P_2$ and $|W|=0$, namely $G\simeq H(n,2,2)$ or $G\simeq H^*(n)$. Some calculation shows that when $4|n$, $H^*(n)$ has one more edge, otherwise they are of the same size.

\begin{figure}[ht]
\begin{center}
\includegraphics{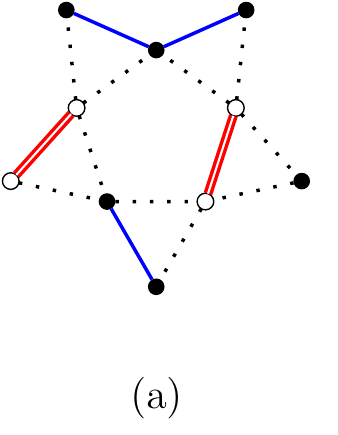}
\hspace{1cm}
\includegraphics{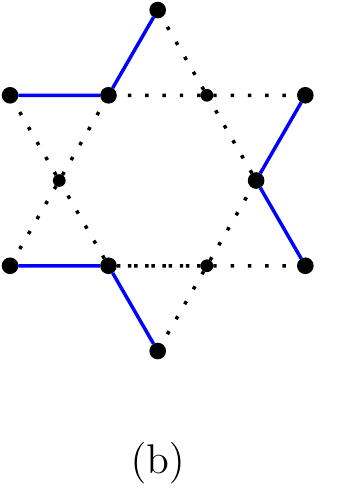}
\hspace{1cm}
\includegraphics{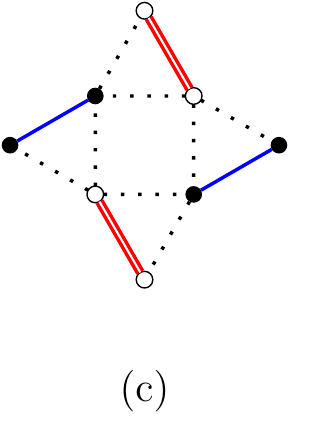}
\end{center}
\caption{(a): $C_k^3\subseteq (P_3\cup \frac{k-3}{2}P_2\cup I_k)\bigotimes (\frac{k-1}{2}P_2)$ when $k=5$; (b): $\frac{k}{2}P_3\in\mathcal{M}(C_k^3)$ when $k=6$; (c): $C_k^3\subseteq \frac{k}{2}P_2\bigotimes \frac{k}{2}P_2$ when $k=4$.}
\label{fig-4}
\end{figure}

If $k\ge 5$, then $t\ge 2$. If $H_1=H_2=P_2$, then since $C_k^3\subseteq (P_3\cup \frac{k-3}{2}P_2\cup I_k)\bigotimes (\frac{k-1}{2}P_2)$ (see Figure~\ref{fig-4}(a)), $W$ is empty, otherwise $P_3\cup \frac{k-3}{2}P_2\cup I_k\subseteq W\cup C_1$. Thus $G\subseteq H^*(n)$, however this implies $e(G)\le e(H^*(n))\sim n^2/4+n/2< H(n,2,t+1)$, a contradiction. Thus $H_1$ and $H_2$ can not both be $P_2$. Suppose $H_1=P_2$ and $H_2=P_1$. Define $W$ to be the set of exceptional vertices that have neighbors in both $A_1$ and $A_2$. Then $|W|<t=\frac{k-1}{2}$, since otherwise $P_{k+1}\subseteq W\cup A_1$ and $P_{k+1}\in\mathcal{M}(C_k^3)$. Then
$$e(G)=\frac{n^2}{4}+\frac{t-1}{2}n+\frac{n}{4}+O(1)<\frac{n^2}{4}+\frac{t}{2}n+O(1)=e(H(n,2,t+1)),$$
a contradiction. Thus $H_1=H_2=P_1$.

When $k$ is even, since $\frac{k}{2}P_3\in\mathcal{M}(C_k^3)$ (see Figure~\ref{fig-4}(b)), $H_i$ is $P_3$-free, $i=1,2$.
Also $C_k^3\subseteq \frac{k}{2}P_2\bigotimes \frac{k}{2}P_2$ (see Figure~\ref{fig-4}(c)), hence at most one $H_i$ has an edge. W.l.o.g. suppose $H_1=P_2$ and $H_2=P_1$. Define $W$ as before and similarly $|W|<t=\frac{k}{2}-1$, which implies $e(G)<e(H(n,2,t+1)),$ a contradiction. Thus $H_1=H_2=P_1$.
\qed
\medskip

For blow-ups of paths, notice that no matter what parity $k$ is, $P_{k+1}^{p+1}\subseteq \lceil \frac{k}{2}\rceil P_2\bigotimes \lceil \frac{k}{2}\rceil P_2$ and $P_{k+1}^{p+1}\subseteq \lceil \frac{k}{2}\rceil P_3\bigotimes I_k$. With these two observations, the same argument works. 
\medskip

\noindent\textbf{Proof of Theorem~\ref{cli-trees} (ii)}
The proof is quite similar to (i), except this time we make use of a forest in the decomposition family obtained by splitting a vertex of degree 2 in $A$. We include here only a sketch of the proof: proceed as in the proof of Theorem~\ref{cli-trees} (i), define $W,A_i,B_i,C_i$ in the same way. Note that still $|W|=a-1$. Since otherwise either the extremal graph $G$ contains a forbidden graph (because some graph in the decomposition family shows up in $G[W\cup C_i]$ for some $i$) or it has fewer edges than $H'(n,p,a)$. It suffices to show that: 

(a) every $G[C_i]$ can have at most one edge, and 

(b) at most one $G[C_i]$ can have one such edge. 

Let $F=T_1\cup T_2\in \mathcal{M}(T^{p+1})$ be the forest obtained by splitting some $z\in A$ with $d(z)=2$. Let $z_1$ and $z_2$ be the two leaves corresponding to $z$ after splitting it and define for $i=1,2$, $T_i'=T_i-z_i$. 

We may assume $a=2$, namely $W$ is non-empty. Since otherwise this tree is a $P_3$, then for (a) two edges in some $G[C_i]$ form a linear forest of size two which is in $\mathcal{M}(P_3^{p+1})$, a contradiction; for (b), if $G[C_i]$ and $G[C_j]$, $i\neq j$, both contains an edge, then we are also done since $P_3^{p+1}\subseteq (P_2\cup I_v)\bigotimes (P_2\cup I_v) \bigotimes T_{n'',p-1}$ for sufficiently large $v$ and $n''$.

For (a), suppose some $G[C_i]$ contains two edges $e_1,e_2$. Similar as Claim~\ref{c2} (ii), $e_1,e_2$ each has at least one endpoint adjacent to $W$. If $e_1,e_2$ are disjoint, then notice that one can embed $T_1'\cup T_2'$ in $G[C_i\cup W]$ and get a copy of $F$ by extending $T_1'$ and $T_2'$ using $e_1$ and $e_2$ respectively. This yields a contradiction since $G[C_i\cup W]$ is $\mathcal{M}(T^{p+1})$-free. If $e_1,e_2$ share an endpoint, namely there is a $P_3=\{w,x,y\}$ in $G[C_i]$. It is not hard to see that there is an embedding for $T$ in $G[C_i\cup W]$, in which $A-\{z\}$ is embedded in $W$ and $z$ is embedded to $x$.

For (b), suppose for $i\neq j$, both $G[C_i]$ and $G[C_j]$ contain an edge $e_i$ and $e_j$ respectively. Then, using $e_i$ and $e_j$, one can partition $W=W_1\cup W_2$, s.t. $T_1\subseteq G[C_i\cup W_1]$ and $T_2\subseteq G[C_j\cup W_2]$, which yields a contradiction since $T^{p+1}\subseteq (T_1\cup I_v)\bigotimes (T_2\cup I_v) \bigotimes T_{n'',p-1}$.
\qed

\end{document}